\documentclass[twoside,11pt]{article}

\pdfoutput=1

\usepackage{calc}
\usepackage{amsmath,amssymb,wasysym}
\usepackage[T1]{fontenc}
\usepackage[utf8x]{inputenc} 
\usepackage[french]{babel}
\usepackage{yhmath} 
\usepackage{enumerate}
\usepackage{xfrac}
\usepackage{url}
\usepackage{layout}
\usepackage{fancyhdr}
\usepackage{textcase}
\usepackage[explicit,rm,it,center,tiny]{titlesec}
\usepackage{graphicx}
\usepackage[colorlinks=true,allcolors=black,hyperfootnotes=false]{hyperref}
\usepackage{array}
\usepackage{ragged2e}
\usepackage[norule,splitrule]{footmisc}
\usepackage{setspace}

\usepackage{arabtex}
\usepackage{utf8}
\title{Un lemme d'analyse dont use Ibn al-Haytham en gnomonique, dioptrique et cinématique céleste}
\author{Erwan Penchèvre}
\date{}
\begin{document}
\maketitle
\footnote{This is a draft version of~: Erwan Penchèvre, «~Un lemme d'analyse dont use Ibn al-Haytham en gnomonique, dioptrique et cinématique céleste~», \textit{Arabic sciences and philosophy}, 2019, 29 (1)~: 133-156. The only difference with the published version is the layout and pagination.}\,Roshdi Rashed a remarqué qu'al-\d{H}asan ibn al-\d{H}asan ibn al-Haytham utilise, dans au moins trois écrits, le lemme suivant~: pour tout $c\in\;]0,1[$ la fonction $x\mapsto f(x)=\dfrac{\sin x}{\sin cx}$ est strictement décroissante sur l'intervalle $]0,\pi/2[$. Nous prendrons prétexte de ce lemme pour présenter certaines \emph{caractéristiques} des travaux d'Ibn al-Haytham dans des domaines \textit{a priori} distincts. Nous nous interrogerons aussi sur la nature même de ce lemme~: Ibn al-Haytham fait-il de l'\emph{analyse}~?
    
Par ce terme, <<~analyse~>>, nous désignerons ici, en un sens à la fois élémentaire et moderne, l'art de <<~majorer, minorer, approcher\,\footnote{Jean Dieudonné a écrit du calcul infinitésimal qu'on <<~pourrait le résumer en trois mots~: majorer, minorer, approcher~>>.}~>>. Pour ne pas courir le risque d'englober toutes les mathématiques sous cette définition, il faut donner à cette discipline son objet, la variable, ou la fonction~; mais peut-être se passerait-on d'une définition explicite du concept de fonction pour reconnaître l'usage opératoire consistant à <<~majorer, minorer~>> une fonction\,\footnote{C'est-à-dire qu'il ne s'agit pas seulement d'encadrer un certain nombre entre deux bornes, mais surtout de démontrer des énoncés \emph{universels} de la forme $(\forall x\in I)\; m<f(x)<M$.}. Cette définition est-elle encore un peu trop vague~? Quand le mathématicien conçoit aujourd'hui l'analyse, il pense aussi à un espace fonctionnel défini par une condition de régularité dont l'archétype est la \emph{continuité}, clairement énoncée et définie au dix-neuvième siècle en termes de limite. Cette évolution entraîne les fondements rigoureux du calcul inventé par Leibniz et Newton, l'exploration du continu et du nombre réel à la fin du dix-neuvième siècle, la théorie de la mesure, etc. C'est donc par rapport à tout cela qu'il nous faudra situer les travaux d'Ibn al-Haytham\,\footnote{Nous nous appuyons sur ces trois traités d'Ibn al-Haytham~: \textit{Sur les lignes des heures}, édité et traduit par R.~Rashed dans \textit{Les Mathématiques infinitésimales du \textsc{ix}\textsuperscript{e} au \textsc{xi}\textsuperscript{e} siècle} (Londres~: Al-Furq\=an Islamic Heritage Foundation, 1996-2006), vol.~V, p.~731-801~; \textit{La Configuration des mouvements des sept astres errants}, \textit{ibid.} p.~263-615~; \textit{Traité sur la sphère ardente}, éd. et trad. par R.~Rashed dans \textit{Geometry and dioptrics in classical Islam} (Londres~: Al-Furq\=an Islamic Heritage Foundation, 2005), p.~224-255. Dans la première partie de notre article (sur les trois applications), nous avons abondamment utilisé les commentaires mathématiques de R.~Rashed à ces trois traités, ainsi que (pour \textit{La Configuration des mouvements}) le commentaire de C.~Houzel dans <<~The new astronomy of Ibn al-Haytham~>>, \textit{A.~S.~P.}~19 (2009)~:~1-41.}.

    \section{Trois applications}

    \subsection{Les lignes des heures saisonnières}
Avant d'en venir à l'étude du lemme et de sa démonstration, nous décrirons trois applications qu'en fait Ibn al-Haytham. La première relève de la théorie gnomonique. Dans son traité \textit{Sur les lignes des heures}, Ibn al-Haytham se demande quelle est la \emph{nature} des lignes des heures sur un cadran solaire. Soit un \emph{cadran solaire horizontal} à la latitude $\varphi=30^\circ$~; c'est à peu près la latitude de Bagdad. On a donc un gnomon vertical et on étudiera la projection centrale par rapport à l'extrémité du gnomon, sur un plan horizontal. La question d'Ibn al-Haytham, triviale si l'on s'intéresse aux heures égales, l'est moins pour les \emph{heures saisonnières}\,\footnote{Chaque heure saisonnière, \textit{zam\=aniya}, ou inégale, \textit{mu`awwaja}, est le douzième de la durée du jour s'étendant du lever au coucher du Soleil.}. Notons par exemple $c=\dfrac{1}{12}$ pour désigner <<~la fin de la première heure~>> (du jour).

    Il existait déjà des ouvrages théoriques sur la projection centrale (ainsi chez al-Qu\d{h}{\=\i} et Ibn Sahl\,\footnote{Voir \textit{Geometry and dioptrics}, et aussi P.~Abgrall, <<~Les débuts de la projection stéréographique~: Conception et principes~>>, \textit{A.~S.~P.} 25 (2015)~:~135-166.}), mais le but poursuivi par Ibn al-Haytham dans ce traité est de répondre de manière théorique à une question pratique\,\footnote{Cf. \textit{Math. inf.} vol.~V p.~684 sur la dialectique pratique / théorique.}. Il souhaite justifier un usage adopté par les constructeurs de cadrans solaires, pour qui les lignes des heures sont des \emph{droites}. Quel est le critère de vérité~? Ce n'est certes pas une question purement mathématique. 

    On assimile l'extrémité du gnomon au centre du Monde, le point $E$ sur la figure~1. Sous cette hypothèse, Ṯābit b. Qurra avait démontré que les images des trajectoires diurnes sont des sections coniques dans le plan du cadran\,\footnote{Ṯābit b. Qurra ne met aucunement en question le domaine de validité de son hypothèse, et il se contente d'écrire~: <<~Nous ne faisons pas de différence entre le point $D$ [ici $E$] et l'extrémité du gnomon~>> (R.~Morelon, \textit{Th\=abit ibn Qurra, {\OE}uvres d'astronomie} (Paris~: Les Belles Lettres, 1987), traité~8, p.~119). Bien que la méthodologie d'Ibn al-Haytham pose la question du domaine de validité des approximations, cette hypothèse-là n'est pas mise en doute dans son traité \textit{Sur les lignes des heures}.}. Au cours de l'année, l'image du Soleil à la fin de la première heure décrit une ligne qui coupe ces sections coniques. Cette ligne est-elle droite~? Selon les Anciens, oui, c'est une droite. Selon Ṯābit, non, c'est une courbe\,\footnote{Voir \textit{Th\=abit ibn Qurra, {\OE}uvres}, traité~9.}. Ibrah{\=\i}m b. Sin\=an démontre que c'est une courbe dans certains cas, et il affirme que c'en est toujours une\,\footnote{Voir le fragment du deuxième livre du traité \textit{Sur les instruments des ombres} dans R.~Rashed et H.~Bellosta, \textit{Ibn Sin\=an, logique et géométrie au \textsc{x}\textsuperscript{e} s.}, p.~414-428.}. Ibn al-Haytham va démontrer que c'est toujours une courbe (sauf dans les cas triviaux)~; mais il réfléchit davantage au critère de vérité ayant cours dans la science dont il traite, et il va en conclure que les Anciens avaient raison~! Voilà l'histoire de ce problème, telle que l'explique Ibn al-Haytham lui-même.

    \begin{figure}
      \footnotesize
      \centering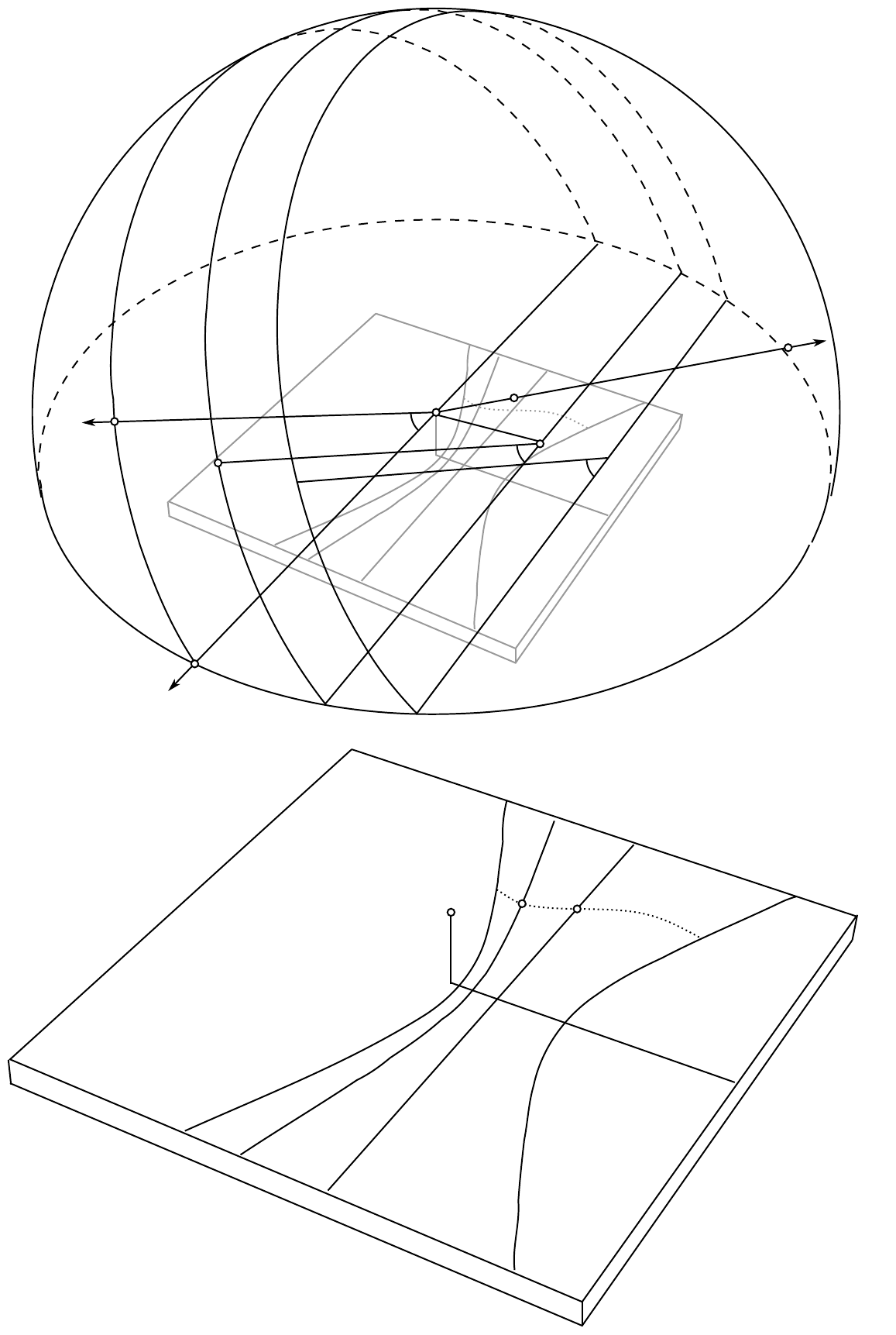

      \normalsize
      \begin{center}
        \textsc{Fig}. 1. Sphère céleste et cadran solaire~; fin de la deuxième heure ($c=2/12$).
      \end{center}
    \end{figure}

\enlargethispage{-\baselineskip}
    Soit donc une trajectoire diurne, le long d'un cercle parallèle à l'équateur céleste. On appellera <<~arc diurne~>> la portion de cette trajectoire située au-dessus de l'horizon, et on notera $y$ sa mesure en radians. Quand le Soleil arrive au point $V$ au douzième de l'arc diurne\,\footnote{Voir fig.~1, où l'on a préféré représenter la fin de la \emph{deuxième} heure ($c=2/12$) pour que la figure soit plus lisible.}, c'est <<~la fin de la première heure~>>. De même, aux équinoxes, à la fin de la première heure, le Soleil arrive au point $H'$ situé au douzième de la moitié visible de l'équateur. On notera $D$ le point Est, et $O$ le centre de la trajectoire diurne considérée. Le plan $H'EV$ coupe le plan de l'équateur et celui de la trajectoire diurne suivant deux droites parallèles $(EH')$ et $(QV)$, où $Q$ est un point situé dans le plan de l'horizon. Ce même plan $H'EV$ coupe le plan du cadran solaire suivant une droite parallèle à la droite $(EQ)$~: la droite joignant les images de $H'$ et $V$ par la projection de centre $E$ sur le plan du cadran\,\footnote{On a noté $p(H')$ et $p(V)$ ces images sur la fig.~1. Quand $y$ varie, la ligne décrite par le point $p(V)$ est la ligne en pointillés. La fig.~1 est exacte, sauf la courbure de la ligne en pointillés qu'on a exagérée pour la rendre sensible.}. \`A $c$ fixé, quand on fait varier la trajectoire diurne, $y$ varie, et Ibn al-Haytham entend donc répondre à la question suivante~: la direction $(EQ)$ reste-t-elle constante~? Si oui, alors la ligne de <<~la première heure~>> sera une droite dans le plan du cadran solaire. On raisonne de même pour toutes les valeurs de $c<1/2$. Les lignes des heures pour $c>1/2$ sont symétriques des précédentes, par rapport au méridien du lieu.

    \begin{figure}
      \footnotesize
      \centering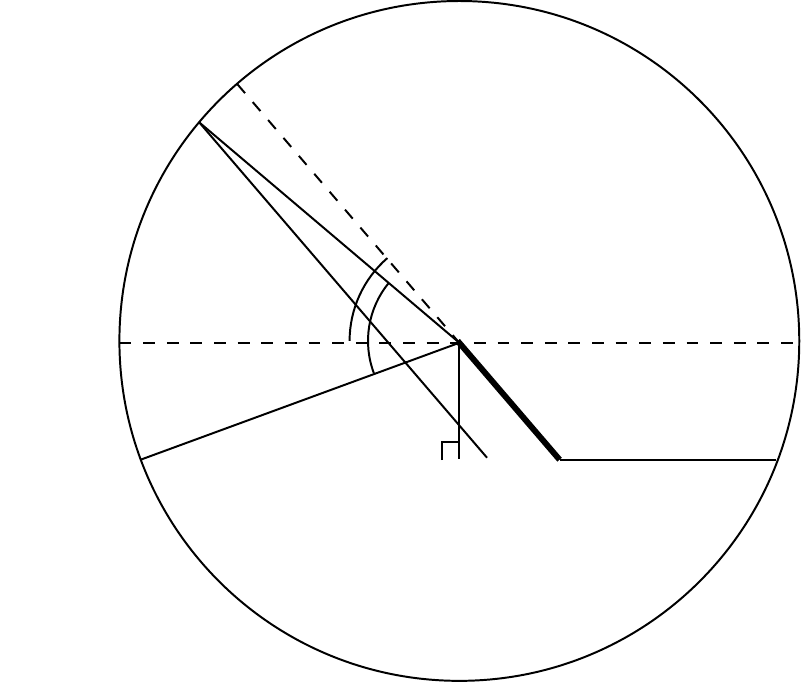

      \normalsize
      \begin{center}
        \textsc{Fig}. 2. Le plan de la trajectoire diurne.
      \end{center}
    \end{figure}
        
    Pour étudier la position de $Q$, Ibn al-Haytham étudie les rapports des segments découpés parallèlement à trois axes $(ED)$, $(EO)$ et $(EH')$. Soit $EO$, $OS$ et $SQ$ ces trois segments (voir fig.~2). Une droite passant par l'origine $E$ est caractérisée par le fait que les rapports des coordonnées affines d'un point sont constants le long de cette droite\,\footnote{Nous ne pensons pas que la notion de coordonnées affines soit totalement étrangère à la pensée d'Ibn al-Haytham. Dans son commentaire mathématique, R. Rashed, peut-être plus prudent, explique ce raisonnement en termes d'homothéties et de triangles semblables (voir \textit{Math. inf.} vol.~V p.~692-729).}. En remarquant que
    $$c\pi+\frac{y-\pi}{2}-cy=\frac{y-\pi}{2}(1-2c),$$
et en s'aidant de la fig.~2, on vérifie facilement que~:
$$\frac{SQ}{OS}=\frac{\sin\left(\dfrac{y-\pi}{2}(1-2c)\right)}{\sin\dfrac{y-\pi}{2}},$$
où $0<1-2c<1$ si $c<1/2$.
  Or notre lemme affirme que ce rapport dépend de manière \emph{strictement croissante} de $y$. Quand le Soleil s'éloigne de l'équateur vers le nord, $y$ croît, $\dfrac{SQ}{OS}$ croît, et la direction $(EQ)$ tend à se rapprocher du méridien~; donc la ligne de la première heure n'est pas droite. Des considérations de symétrie montrent que l'image de $H'$ est un point d'inflexion autour duquel serpente la courbe comme dessinée en pointillés.

  $\dfrac{SQ}{OS}$ atteint son maximum au solstice d'été quand, à Bagdad~:
  $$\dfrac{y-\pi}{2}=\dfrac{y_{\max}-\pi}{2}\simeq 15^\circ.$$
  D'autre part, la fonction $x\mapsto\dfrac{\sin x}{x}$ est décroissante sur l'intervalle $]0,\pi/2[$ (cas limite du lemme, déjà connu de Ptolémée\,\footnote{Ptolémée, \textit{The Almagest}, trad. R.~Catesby Taliaferro, coll. <<~Great books of the Western world~>>, vol.~16, \textit{Encyclop{\ae}dia Britannica}, 1952, I.10 p.~19. Comme nous l'a fait remarqué Marc Dehon, Archimède avait déjà énoncé ce cas limite dans \textit{L'Arénaire} -- voir \textit{Archimède}, texte établi et trad. par C.~Mugler (Paris~: Les Belles Lettres, 1970-1972), tome~II, p.~143.})~:
      $$\frac{\sin\dfrac{y-\pi}{2}}{\dfrac{y-\pi}{2}}\leq\frac{\sin\left(\dfrac{y-\pi}{2}(1-2c)\right)}{\dfrac{y-\pi}{2}(1-2c)},$$
    d'où l'encadrement suivant~:
    $$1-2c \leq \dfrac{SQ}{OS} \leq \dfrac{\sin\left(\dfrac{y_{\max}-\pi}{2}(1-2c)\right)}{\sin\dfrac{y_{\max}-\pi}{2}}.$$

    \begin{figure}
      \footnotesize
      \centering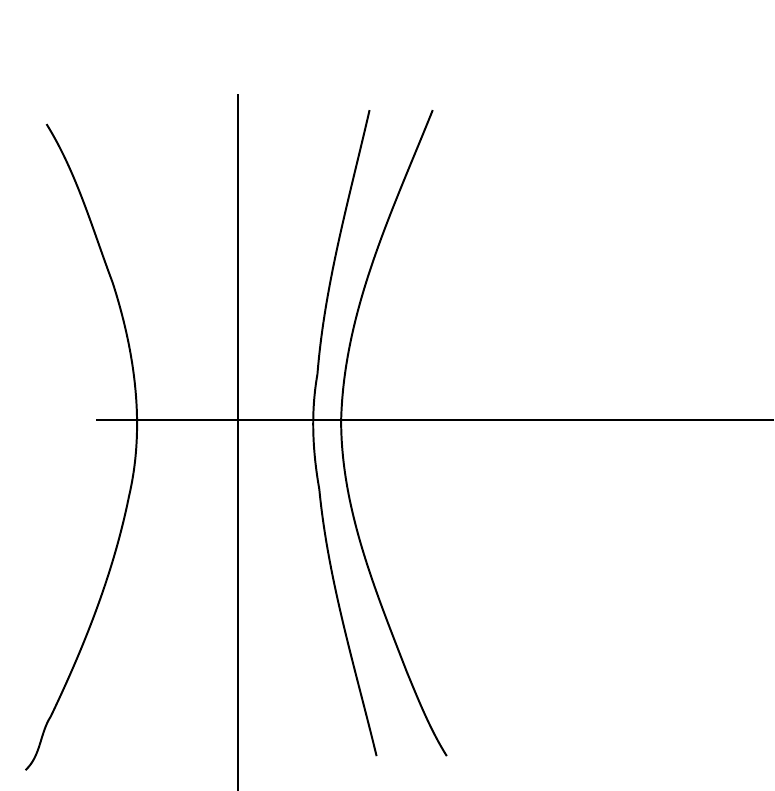

      \normalsize
      \begin{center}
        \textsc{Fig}. 3. Le plan du cadran (on a beaucoup exagéré la courbure en pointillés).
      \end{center}
    \end{figure}
    
    Ibn al-Haytham ne s'arrête pas à cet encadrement~; il rapporte aussi la droite image $(IQ)$ à deux directions prises dans le plan du cadran (voir à présent fig.~3 pour les notations\,\footnote{Pour garder les notations d'Ibn al-Haytham, nous avons noté $Q=p(V)$ dans le plan du cadran fig.~3, bien que la lettre $Q$ ait été déjà utilisée ci-dessus pour désigner un autre point dans le plan de l'horizon~; la lettre $O$ désigne aussi un nouveau point~; et nous avons choisi $V$ sur une trajectoire diurne située du côté du Capricorne, contrairement aux fig.~1 et 2.}), et il majore le rapport $\dfrac{TQ}{IT}$, pour $c=\dfrac{k}{12}$, $0<k<12$, à la latitude $\varphi\simeq 30^\circ$ de Bagdad. Il montre que\,\footnote{Cette valeur n'est autre que la largeur de l'encadrement précédent, multipliée par la valeur $\left(\dfrac{OS}{EQ}\right)_{\max}$ de $\dfrac{OS}{EQ}$ quand $y=y_{\max}$. On passe les détails.}~:
    $$\dfrac{TQ}{IT}<\dfrac{LO}{IL}=\dfrac{1}{174}.$$
    Ceci lui permet de comparer $TQ$ à un \emph{étalon physique}. Pour un gnomon typique d'une longueur de trois doigts, il calcule que $TQ$ est inférieur à $\dfrac{1}{30}\times 3\text{ doigts}=\dfrac{3}{5}$ de la longueur d'un grain d'orge. Cette grandeur étant négligeable par rapport à la dimension caractéristique du cadran solaire, il en conclut que la ligne des heures \emph{n'est pas} une droite mathématique \emph{pour l'imagination}, mais qu'elle \emph{est} droite \emph{pour la sensation}.

    Autrement dit, fig.~3, toutes les droites comprises entre $IL$ et $IO$ sont une même droite pour la sensation. La droite pour la sensation a une certaine \emph{épaisseur}, et cette épaisseur peut faire l'objet d'une \emph{mesure}. Voilà l'objet d'une théorie physique, et son procédé est l'analyse~-- au sens de <<~majorer, minorer, approcher~>>\,\footnote{R.~Rashed y voit même les débuts d'une <<~théorie des erreurs~>> où l'on recherche des encadrements permettant de choisir les paramètres à négliger en ordre de grandeur (\textit{Math. inf.} vol.~V p.~692).}.

    \begin{figure}
      \footnotesize\centering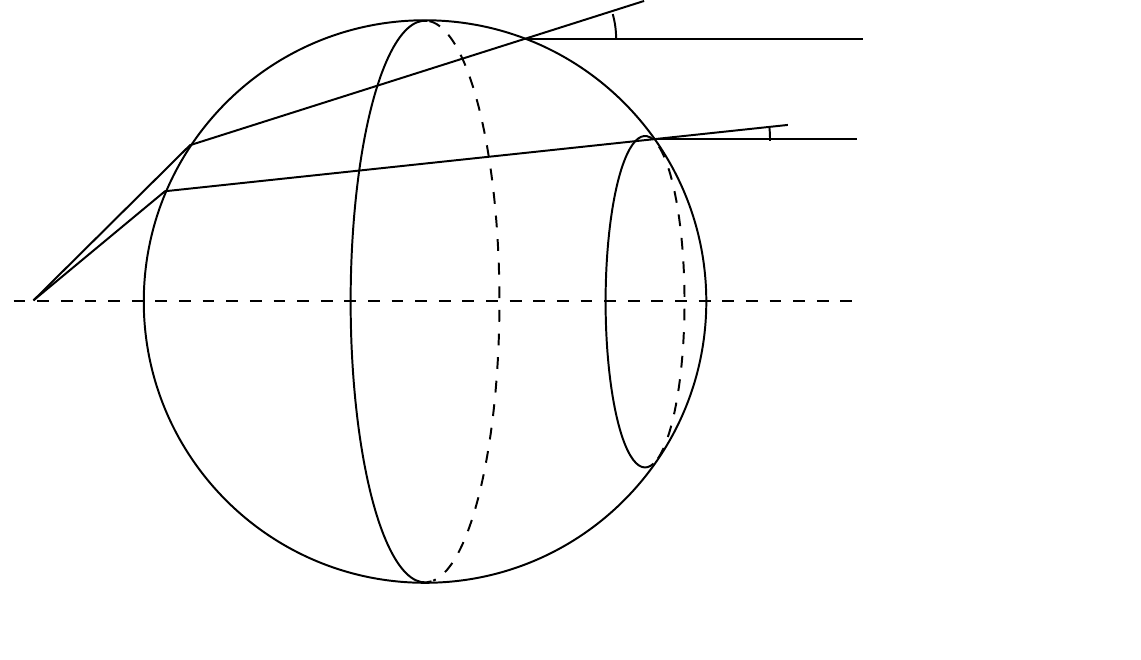

      \normalsize
      \begin{center}
        \textsc{Fig}. 4. Double réfraction d'un faisceau parallèle par un dioptre sphérique.
      \end{center}
    \end{figure}

    \subsection{Double réfraction d'un faisceau parallèle par un dioptre sphérique}
Ibn al-Haytham fait une autre application du lemme, dans son étude de la double réfraction d'un faisceau de rayons parallèles par un dioptre sphérique~: dans le traité \textit{Sur la sphère ardente}. Il n'avait pas abordé ce cas dans son \textit{Optique}\,\footnote{Voir \textit{Geometry and dioptrics} p.~153-159 pour un résumé des cas de réfraction par un dioptre sphérique étudiés dans l'\textit{Optique} d'Ibn al-Haytham.}. Ibn al-Haytham pose deux questions (voir fig.~4)~: 1) où est le foyer ardent~? 2) deux rayons d'incidences distinctes, parallèles à un axe du dioptre, peuvent-ils être réfractés vers un même point de l'axe~? On verra que ces deux questions sont liées. On commencera par rappeler quelques lois de la réfraction qu'Ibn al-Haytham tire de son \textit{Optique}~:
    
    \begin{enumerate}[a)]
    \item L'angle de déviation $d$ et l'angle de réfraction $r$ sont deux fonctions croissantes de l'angle d'incidence $i$. 
    \item Si $r$ est l'angle de réfraction associé à l'angle d'incidence $i$ à l'entrée du dioptre, alors, à la sortie du dioptre, sous une incidence $r$, l'angle de réfraction sera égale à $i$ (<<~retour inverse~>>)\,\footnote{Les points a) et b) impliquent que, si deux rayons parallèles, après double réfraction par le dioptre sphérique, pouvaient converger vers un même point de l'axe, alors ils ne pourraient pas se croiser à l'intérieur du dioptre.}.
    \item Dans l'air, à l'entrée d'un dioptre en verre, $0<i-2d< 2d$. On sait aujourd'hui que la seconde inégalité est vraie, mais la première n'est vraie que jusqu'à une certaine incidence~: seulement si $i/2<\arccos(3/4)$. C'est donc seulement en négligeant les rares incidences supérieures que les résultats ci-dessous seront valables.
    \item La fonction $i\mapsto i/d$ est décroissante\,\footnote{C'est une conséquence de la loi de Snell-Descartes, mais cette proposition avait déjà été énoncée par Ptolémée. \`A cette occasion, rappelons que la loi de Snell-Descartes était connue d'Ibn Sahl au \textsc{x}\textsuperscript{e} siècle~; Ibn al-Haytham n'a pas adoptée cette loi découverte par son prédécesseur, peut-être parce qu'elle n'était pas, chez Ibn Sahl, fondée sur l'observation et l'expérimentation, cf. \textit{Geometry and dioptrics} p.~177-182 et 1039-1045. On peut remarquer que la démonstration de d) à partir de la loi de Snell-Descartes repose sur le lemme dont il est ici question~: si $i_1<i_2$, alors
      $$\frac{\sin i_1}{\sin(1-d_1/i_1)i_1}=\frac{\sin i_1}{\sin r_1}=\frac{\sin i_2}{\sin r_2}=\frac{\sin i_2}{\sin(1-d_2/i_2)i_2}\stackrel{\text{(lemme)}}{<}\frac{\sin i_1}{\sin(1-d_2/i_2)i_1},$$
    d'où $d_1/i_1<d_2/i_2$.}.
    \end{enumerate}

    \begin{figure}
      \footnotesize\centering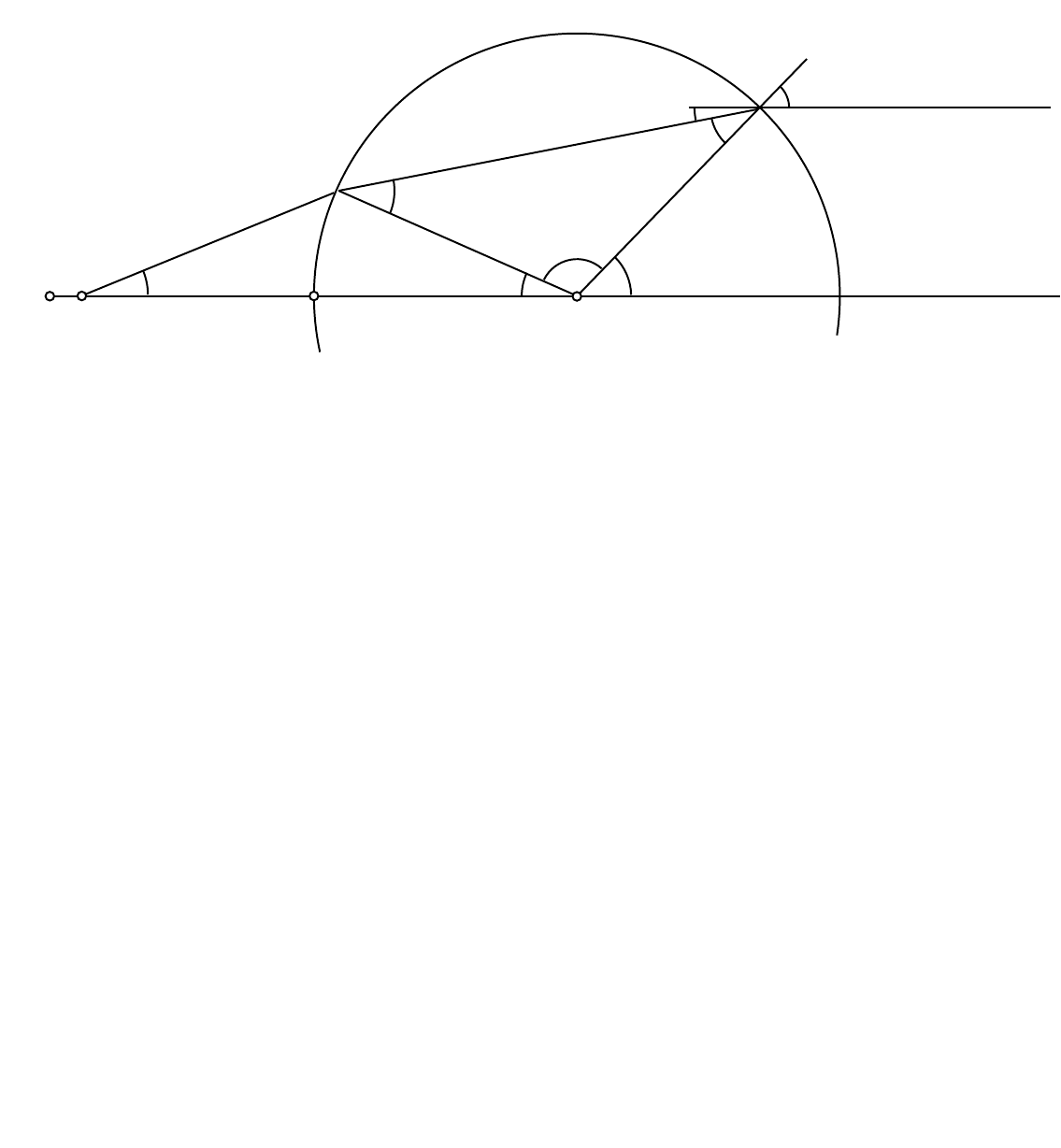

      \normalsize\begin{center}
        \textsc{Fig}. 5. Démonstration par l'absurde.
      \end{center}
    \end{figure}
    
    Ibn al-Haytham montre que la réponse à la question 2) est négative\,\footnote{C'est la proposition 3 du traité, voir \textit{Geometry and dioptrics} p.~232.}, en démontrant par l'absurde que la fig.~5~(ii) est impossible, si $i_1<i_2$ sont les angles d'incidences de deux rayons parallèles à l'axe en entrée du dioptre\,\footnote{Pour comprendre la fig.~5~(ii), voir d'abord la fig.~5~(i) qui illustre la double réfraction d'un unique rayon d'incidence $i$ en entrée du dioptre. On remarque que $r+d=i\Rightarrow 2r-i=i-2d$, et la loi c) implique que le point où le rayon sortant rencontre l'axe est à une distance du dioptre inférieure au rayon de la sphère.}. Si en effet elle était possible, alors on aurait~:
    $$\frac{\sin(2d_2+2d_1)}{\sin(2d_2-2d_1)}=\frac{TR}{RO}>\frac{TW}{WO}=\frac{\sin(i_2-2d_2+i_1-2d_1)}{\sin((i_2-2d_2)-(i_1-2d_1))}$$
    Mais d) implique que~:
    $$\frac{i_2-2d_2}{2d_2}<\frac{i_1-2d_1}{2d_1},$$
    d'où\,\footnote{$\dfrac{a}{b}<\dfrac{c}{d}\Rightarrow\dfrac{b+d}{b-d}<\dfrac{a+c}{a-c}$.}
    $$\frac{2d_2+2d_1}{2d_2-2d_1}<\frac{i_2-2d_2+i_1-2d_1}{(i_2-2d_2)-(i_1-2d_1)}.$$
    S'il y avait égalité, on aurait (à cause du lemme)~:
    $$\frac{\sin(2d_2+2d_1)}{\sin(2d_2-2d_1)}<\frac{\sin(i_2-2d_2+i_1-2d_1)}{\sin((i_2+2d_2)-(i_1-2d_1))},$$
    mais il n'y a pas égalité, le membre droit est supérieur, donc c'est encore pire\dots\ Absurde.

    Pour un faisceau de rayon parallèles, au vu de la réponse à la question 2), chaque point de l'axe ne reçoit donc qu'un unique rayon du faisceau. Si le \emph{foyer ardent} est un point en lequel convergent un plus grand nombre de rayons à la sortie du dioptre, il est donc clair qu'un tel foyer n'existe pas.

    \enlargethispage*{\baselineskip}
    Mais à nouveau, puisque la dioptrique n'est pas, pour Ibn al-Haytham, une mathématique pure, il faut prendre garde au critère de vérité. Le \emph{point mathématique} n'est peut-être pas adéquat à la réalité physique. Comme le dit Ibn al-Haytham, un point mathématique n'a pas de volume, il ne peut donc pas s'échauffer\,\footnote{Voir \textit{Geometry and dioptrics} p.~250-251~:
      \begin{quote}
        \begin{RLtext}
          \setcode{utf8}
        
       ومن أجل أن هذا الموضع ذو مقدار, صارت فيه حرارة. ولو كانت نقطة متوهمة, لما حصل فيها حرارة.
        
     \end{RLtext}
     \end{quote}
     Ibn al-Haytham ne saurait confondre la droite mathématique et la lumière~-- phénomène physique. Si la droite est conçue comme élément idéal, le phénomène lumineux n'a d'effet qu'en intégrant ces éléments idéaux en un faisceau de rayons. C'est un trait de la science d'Ibn al-Haytham que la lumière est pour lui un objet physique \emph{distinct} de l'idéalité mathématique mais pouvant être étudié mathématiquement. Cf. l'article <<~Des raisons de la lumière aux lumières de la raison~>>, contribution au colloque <<~The Islamic golden age of science for today's knowledge-base society~: The Ibn al-Haytham example~>>, 14-15 sept. 2015, UNESCO, Paris, où Michel Paty compare ce passage d'une théorie mathématique à une théorie physique mathématisée avec la prise en compte par Einstein du caractère physique de la lumière en relativité restreinte et en quantique.}. Inutile de compter le nombre de rayons convergeant exactement en un point --~mathématique~-- de l'axe~; au contraire, il faut mesurer de petits faisceaux constitués chacun d'une infinité --~continue~-- de rayons, et comparer entre eux de tels faisceaux atteignant de petites zones --~non strictement ponctuelles~-- de l'axe. Mais comment comparer deux faisceaux de rayons, en l'absence d'une théorie de la mesure, d'une notion d'angle solide, ou plutôt ici, de flux à travers une surface~?

    Peut-être faut-il ici bien peser le geste encore implicite dans la solution proposée par notre savant. Ibn al-Haytham se contente d'isoler deux faisceaux de rayons incidents parallèles à l'axe en entrée, et d'étudier les points où ils rencontrent l'axe en sortie du dioptre~: les rayons d'incidence supérieure à $50^\circ$, et ceux d'incidence inférieure à $40^\circ$. Il semble croire, sans justification rigoureuse, que ces deux faisceaux convoient la même quantité de chaleur. Et en effet, sur un plan perpendiculaire à l'axe, avant d'entrer dans le dioptre, ces deux faisceaux interceptent deux aires égales~: cf. les deux aires hachurées sur la fig.~6~(i). Ibn al-Haytham l'avait sûrement calculé.

    \begin{figure}
      \footnotesize\centering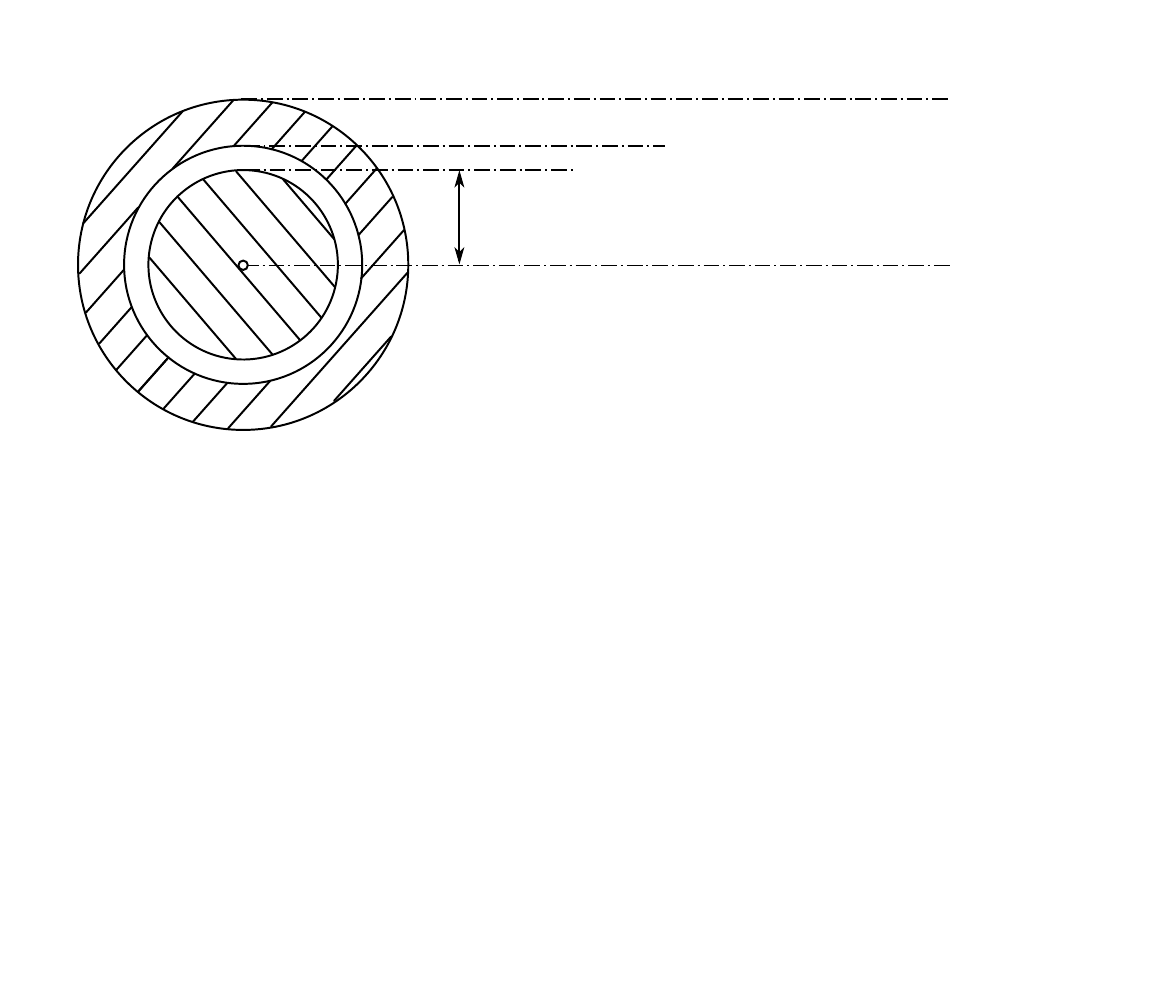

      \normalsize\begin{center}
        \textsc{Fig}. 6. Comment comparer deux faisceaux.
      \end{center}
    \end{figure}
    
    Il faut à présent étudier les rayons sortants. Notons $CD=1$ le rayon du dioptre, $V$ le point de l'axe tel que $VC=CD$, et $S$ le milieu de $VC$ (voir fig.~8 pour les notations). Il existe deux points $N$ et $N'$ de l'axe tels que les rayons d'incidence $\geq 50^\circ$ sont réfractés vers $N'C$, et les rayons d'incidence $\leq 40^\circ$ vers $VN$. Ptolémée avait montré que pour $i=40^\circ$, on a $d=15^\circ$, $2d=30^\circ$, $i-2d=10^\circ$, et~:
    $$DX=\cos(10^\circ),$$
    $$KX=\sin(10^\circ),$$
    $$\frac{KX}{NX}=\tan(30^\circ)=\frac{\sqrt{3}}{3},$$
    $$DN=DX+NX=\cos10^\circ+\sqrt{3}\sin10^\circ\simeq1,3<1,5.$$
    Ibn al-Haytham en conclut que <<~les rayons réfractés vers $SC$ sont bien plus nombreux que les rayons réfractés vers $SV$~>>\,\footnote{Voir \textit{Geometry and dioptrics}, p.~250.}, et en effet, rien que sur $N'C$, il en arrive autant que sur $NV$, donc plus que sur $SV$. C'est une réponse quantitative à la question 1)~: le foyer ardent est à une distance du dioptre inférieure au quart de son diamètre.

    \enlargethispage*{\baselineskip}
    \subsection{Cinématique céleste}
    Il nous reste à mentionner une troisième application du lemme. Ce lemme est au fondement du long traité d'Ibn al-Haytham sur \textit{La configuration des mouvements des sept astres errants}. Il serait impossible de décrire ici tout le contenu de ce traité de manière aussi détaillée que nous l'avons fait pour les deux textes précédents. Dans ce traité, Ibn al-Haytham entreprend l'étude mathématique des trajectoires des astres dans un référentiel attaché à l'observateur. Traditionnellement, le référentiel priviligié pour les théories planétaires était un référentiel sidéral (chez Ptolémée) ou bien tropique\,\footnote{Tropique, peut-être depuis les Ban\=u M\=us\=a qui <<~repoussent l'écliptique au-delà de la sphère des fixes~>> comme l'explique R.~Morelon dans son commentaire du \textit{Livre sur l'année solaire}, dans \textit{Th\=abit ibn Qurra, {\OE}uvres}, p.~\textsc{lxvi}.}. Ce qui intéresse Ibn al-Haytham, c'est le phénomène. Comme en optique, il faut penser le phénomène non plus comme un objet mathématique idéal, mais en quelque sorte en \emph{intégrant} de tels objets idéaux. La trajectoire d'un astre est pensée alors comme une succession de déplacements d'un point à un autre~; les accroissements finis des coordonnées seront soumis au calcul et à l'étude quantitative. C'est très différent du concept traditionnel de mouvement céleste, combinaison de quelques mouvements simples --~souvent une composée d'un nombre fini de rotations uniformes~-- dont l'étude mathématique ne décrivait jamais qu'\emph{une seule position} de l'astre, à \emph{un seul instant}, sauf pour des questions particulières\,\footnote{On pense ici à la question des coascensions des arcs de l'écliptique, \textit{ma\d{t}\=ali` al-bur\=uj}~; Neugebauer y voyait l'une des motivations majeures de la géométrie sphérique hellénistique (Théodose, Ménélaüs), cf. \textit{The Exact sciences in Antiquity} (N.~Y.~: Dover, 1969), p.161.}.

    Cette nouvelle théorie va, comme dans les deux textes ci-dessus, se heurter au paradoxe~: le mathématicien, en manipulant les objets idéaux dans l'imagination, va prédire des phénomènes dont les dimensions caractéristiques échappent au domaine de l'expérience ou de l'observation. L'analyse devra alors déterminer, par des encadrements, l'ordre de grandeur de ces phénomènes, et elle devra les comparer à des étalons physiques.

    \textit{La configuration des mouvements} --~ou ce qui nous en est parvenu~-- pose un tel problème précis. Ce problème oriente toute la succession des lemmes et des propositions, et il se décline en les deux variantes suivantes~:
    \begin{enumerate}[1)]
      \item Lorsque la trajectoire diurne d'un astre frôle l'horizon de l'observateur, l'astre va se lever immédiatement après s'être couché. Est-il possible que l'astre se couche à l'est du méridien (ou bien qu'il se lève à l'ouest)\,\footnote{Deux siècles après Ibn al-Haytham, Na\d{s}{\=\i}r al-D{\=\i}n al-\d{T}\=us{\=\i} mentionne encore ce problème au rang des <<~questions étranges et difficiles~>> (voir F.~J. Ragep, \textit{Na\d{s}{\=\i}r al-D{\=\i}n al-\d{T}\=u\d{s}{\=\i}'s Memoir on astronomy} (N.~Y.~: Springer-Verlag, 1993), III.5 [7] p.~278). Ibn al-Š\=a\d{t}ir le mentionne aussi, voir E.~Penchèvre, <<~La \textit{Nih\=aya al-s\=ul f{\=\i} ta\d{s}\d{h}{\=\i}\d{h} al-'u\d{s}\=ul} d'Ibn al-Š\=a\d{t}ir~: \'Edition, traduction et commentaire~>>, à paraître, fin du chapitre II.5 p.~348-349.}~?
      \item Chaque astre ayant un lever et un coucher culmine près du méridien une fois chaque jour. Est-il possible que la hauteur maximale à l'horizon soit en fait atteinte un peu à l'est du méridien (ou bien un peu à l'ouest)~?
    \end{enumerate}

    \begin{figure}
      \footnotesize\centering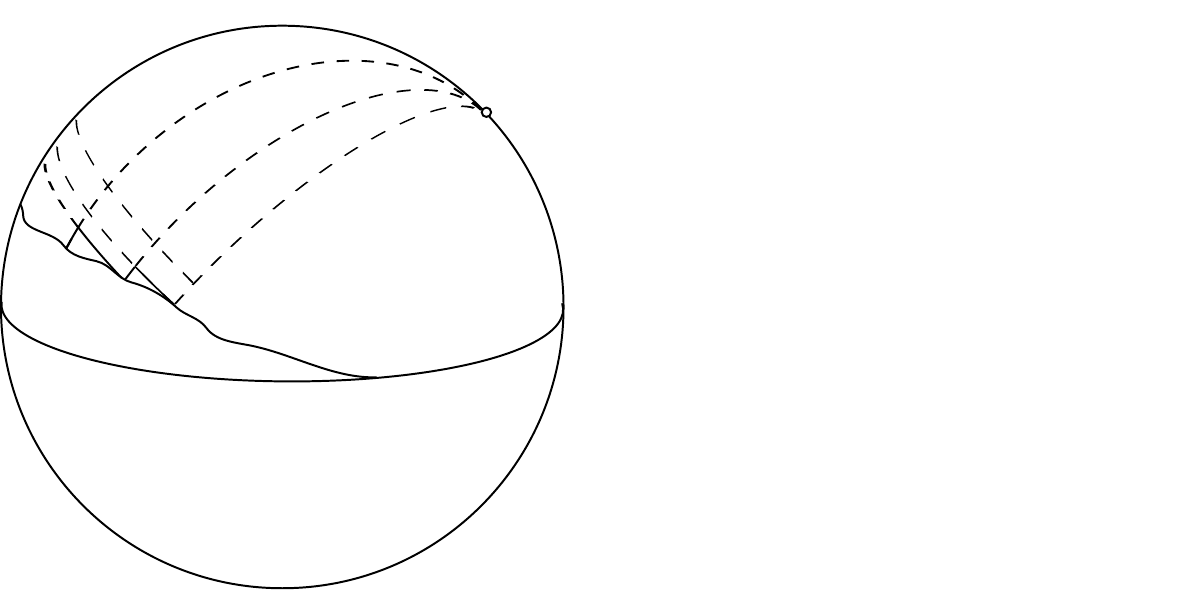

      \normalsize\begin{center}
        \textsc{Fig}. 7. Trajectoires dans un référentiel attaché à l'observateur.
      \end{center}
    \end{figure}
    
    Pour répondre à ces questions de hauteurs sur l'horizon, Ibn al-Haytham choisit un système de coordonnées sphériques <<~équatoriales~>>, mais dans un référentiel attaché à l'observateur~: \emph{déclinaison} et \emph{temps requis}\,\footnote{\`A un instant donné, le \emph{temps requis} est mesuré comme une \emph{ascension dans la sphère droite}~; mais, à la différence du temps requis, l'ascension droite repère la position des astres dans un référentiel \emph{en mouvement} par rapport à l'observateur (le mouvement diurne).}. Le pôle Nord est supposé immobile par rapport à l'observateur. En fait, Ibn al-Haytham ne précise pas l'origine de la coordonnée \emph{temps requis} car il ne considère jamais que des \emph{accroissements} en temps requis. Cela permet au savant d'imaginer des figures qui devaient ressembler à la fig.~7~(i), où $\Delta$ désigne un accroissement en déclinaison et $\delta$ un accroissement en temps requis.

    \`A cause du fait que le mouvement diurne domine très nettement (même pour la Lune) les mouvements propres des astres par rapport à l'écliptique, pour l'observateur, la trajectoire apparente d'un astre progresse toujours d'est en ouest, c'est-à-dire qu'elle est paramétrisable selon le temps requis. L'étude des trajectoires revient donc à étudier des taux d'accroissements finis $\Delta/\delta$. Le c{\oe}ur du traité consiste à déduire des modèles planétaires une minoration de $\Delta/\delta$ le long d'arcs de trajectoire évitant les extrémités nord et sud~: les extrémités nord et sud d'une trajectoire sont en effet des extréma en déclinaison où $\Delta/\delta$ tend à s'annuler. Ce n'est pas le moindre mérite d'Ibn al-Haytham que d'avoir su cerner l'enjeu et la possibilité de telles minorations.

    \enlargethispage*{\baselineskip}
    Nous nous contenterons d'expliquer ici la proposition 28 du traité d'Ibn al-Haytham\,\footnote{Voir \textit{Math. inf.} vol.~V p.~524-537.}. Soit un arc de trajectoire, allant vers le sud, le long duquel l'astre se lève à l'est et passe au méridien~; alors il existe un maximum de hauteur strictement avant le passage au méridien.

    Voici l'idée de la démonstration. Supposons qu'on ait une minoration\linebreak ${\Delta/\delta>m}$ le long d'un tel arc de trajectoire. Soit $B$ le point de l'horizon où l'astre se lève, et $D$ le point où il passe au méridien, voir fig.~7~(ii). Ibn al-Haytham construit un point $I$ sur le méridien, tel que
    $$\frac{\wideparen{ID}}{\wideparen{HI}}<m,$$
    où $H$ est le point situé à même hauteur que $D$ et de même déclinaison que $I$. Le point $H$ est alors l'intersection de l'almucantar de hauteur de $D$ et du parallèle passant par $I$. Ici $\wideparen{HI}$ désigne l'angle au centre dans le cercle $HI$, et non la longueur de l'arc\,\footnote{Ibn al-Haytham l'explique ainsi~: <<~Quand je dis “\,le rapport de l'arc\,”, dans deux cercles différents, j'entends le rapport de l'arc (du cercle égal au cercle) semblable à l'arc du premier cercle~; c'est cela que j'entends dans tout ce que je mentionnerai plus tard si j'utilise les rapports des arcs les uns aux autres.~>> (\textit{Math. inf.} vol.~V p.~350-351.)}~; de même pour $\wideparen{ID}$. Que la trajectoire traverse le parallèle passant par $I$, au point $M$. Considérons l'accroissement fini $\dfrac{\Delta}{\delta}$ du point $M$ au point $D$. On a $\Delta=\wideparen{ID}$, et $\dfrac{\wideparen{ID}}{\wideparen{HI}}<m<\dfrac{\Delta}{\delta}$, donc $\delta<\wideparen{HI}$. Ainsi le point $M$ est entre $H$ et $I$, et la trajectoire atteint donc une hauteur strictement supérieure à celle $D$, un peu avant le passage au méridien.

    Pour le Soleil, $\wideparen{ID}$ ne dépasse pas quelques minutes d'arc. Pour la Lune, il peut atteindre environ $1^\circ$, mais rarement\,\footnote{Voir \textit{Math. inf.} vol.~V p.~560.}. Le phénomène physique est donc peu sensible. L'importance de ces majorations tient peut-être aussi au fait mathématique suivant. Ibn al-Haytham, à la fin de la partie du traité ayant survécu, cherche à démontrer l'\emph{unicité} du maximum de hauteur, et il est alors amené à assimiler de petits triangles sphériques à des triangles plans~; mais nous ne pouvons examiner ici cette partie du texte.
    
    \section{Deux démonstrations}

    \subsection{Démonstration locale}
Venons-en à la démonstration du lemme proposée par Ibn al-Haytham dans son traité \textit{Sur les lignes des heures}. Dans ce traité, c'est le lemme 3. Nous admettrons le cas limite déjà connu de Ptolémée~: la fonction $x\mapsto f(x)=\dfrac{\sin x}{x}$ est strictement décroissante sur l'intervalle $]0,\pi/2[$. On remarque que ceci vaut aussi pour les cordes, car~:
    $$f(x)=\frac{\sin x}{x}=\frac{\frac{1}{2}\text{crd}\,2x}{\frac{1}{2}2x}=\frac{\text{crd}\,2x}{2x}.$$

    Avant tout, Ibn al-Haytham énonce un premier lemme (lemme 1) qui revient à formuler l'encadrement suivant, pour $0<\alpha<\beta<\pi/2$~:
    $$\frac{\beta}{\alpha}\in\left[\frac{\sin\beta}{\sin\alpha},\frac{\sin\frac{\alpha+\beta}{2}+\sin\frac{\alpha-\beta}{2}}{\sin\frac{\alpha+\beta}{2}-\sin\frac{\alpha-\beta}{2}}\right].$$
    \emph{Démonstration~:} La minoration n'est autre que le cas limite connu de Ptolémée. Quant à la majoration, on l'obtient aussi à partir de ce cas limite, en posant ${u=\dfrac{\beta-\alpha}{2}}$,\enskip $v=\dfrac{\beta+\alpha}{2}$. On a
    $$\frac{\beta}{\alpha}=\frac{2}{1-u/v}-1.$$
    Or $v>u$, et $\dfrac{v}{u}>\dfrac{\sin v}{\sin u}$, donc~:
    $$\frac{\beta}{\alpha}<\frac{2}{1-\frac{\sin u}{\sin v}}-1.$$
    D'où la majoration souhaitée.
    
    Quant au lemme 3, Ibn al-Haytham commence par en démontrer un cas particulier facile (c'est le lemme 2 du traité)~: le cas $c=1/2$. Dans ce cas $f(x)=\dfrac{\sin x}{\sin(x/2)}=2\cos(x/2)$.
    
    Pour le cas général (lemme 3 du traité), nous tenterons ici de restituer l'analyse\,\footnote{Ici, nous entendons \emph{analyse} au sens de la géométrie euclidienne, c'est-à-dire qu'on <<~lit la démonstration en partant de la fin~>>, etc. Pour ne pas allonger inutilement, on renvoie à la fig.~8 pour toutes les notations, et on passera certains détails sous silence, pour lesquels on renvoie à la synthèse (dans le texte original ou bien dans le commentaire qu'en donne R. Rashed). On remarquera que, en général, comme on emploiera certaines majorations grossières, on ne progressera pas toujours par équivalence logique pour passer d'une inégalité souhaitée à la suivante.} qu'a dû faire Ibn al-Haytham. Soit donc $0<x<y\leq\pi/2$. On veut démontrer que~:
    $$\frac{\sin x}{\sin cx}\stackrel{?}{>}\frac{\sin y}{\sin cy}$$
    c'est-à-dire, quitte à représenter $x$ et $y$ comme des angles au centre d'un même cercle~-- voir fig.~8 pour les notations~:
    $$\frac{\sin\wideparen{BD}}{\sin\wideparen{BE}}\stackrel{?}{>}\frac{\sin\wideparen{BA}}{\sin\wideparen{BC}}$$
    ou encore~:
    $$\frac{DI}{IE}\stackrel{?}{>}\frac{AH}{HC}.$$

    \begin{figure}
      \footnotesize\centering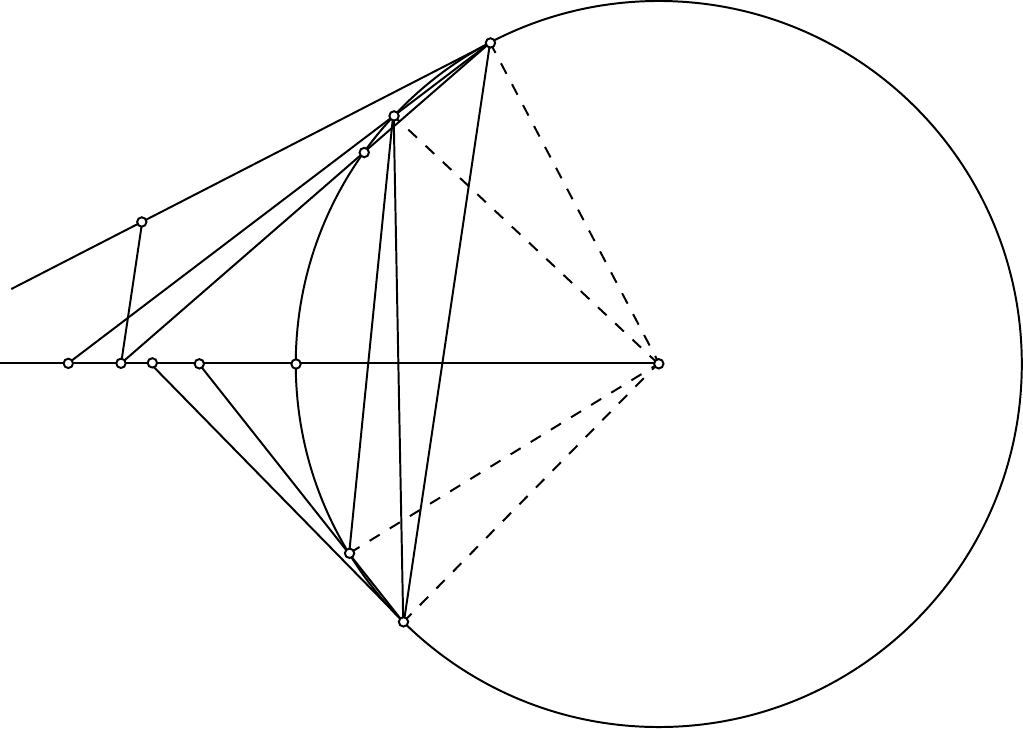
      
      \normalsize\begin{center}
        \textsc{Fig}. 8. Démonstration du lemme dans le traité \textit{Sur les lignes des heures}.
      \end{center}
    \end{figure}

    C'est ici qu'intervient une construction auxiliaire non triviale~: les tangentes au cercle en $A$ et $C$. Heuristiquement, on peut comprendre cette étape si l'on imagine $x$ proche de $y$. Puis Ibn al-Haytham applique deux fois le théorème de Ménélaüs pour se ramener d'une comparaison de rapports le long de transversales au cercle à une comparaison de rapports le long de transversales plus proches des tangentes\dots\ On espère ensuite pouvoir comparer les transversales aux tangentes. Le théorème de Ménélaüs donne~:
    $$\frac{DI}{IE}=\frac{DW}{WC}\cdot\frac{CR}{RE}\text{\enskip et\enskip}\frac{AH}{HC}=\frac{AU}{UD}\cdot\frac{DW}{WC}.$$
    Il reste donc à démontrer que~:
    $$\frac{CR}{RE}\stackrel{?}{>}\frac{AU}{UD}$$
    $$1/\left(1-\frac{CE}{CR}\right)\stackrel{?}{>}1/\left(1-\frac{AD}{AU}\right)$$
    $$\frac{CE}{CR}\stackrel{?}{>}\frac{AD}{AU}$$
    Le cas limite dû à Ptolémée, appliqué aux cordes $CE$ et $AD$, montre que $\dfrac{CE}{\wideparen{CE}}>\dfrac{AD}{\wideparen{AD}}$. Or $(CE/\wideparen{CE})\cdot(\wideparen{CE}/CR)=CE/CR$ et $(AD/\wideparen{AD})\cdot(\wideparen{AD}/AU)=AD/AU$. Il suffit donc de démontrer~:
    $$\frac{\wideparen{CE}}{CR}\stackrel{?}{>}\frac{\wideparen{AD}}{AU}.$$
    Ici, Ibn al-Haytham utilise la majoration grossière $CR<CJ$ (car l'angle $\widehat{CRJ}$ est obtus). On est réduit à montrer que~:
    $$\frac{\wideparen{CE}}{CJ}\stackrel{?}{>}\frac{\wideparen{AD}}{AU}.$$
    Si l'on fixe les points $A$, $B$ et $C$, et qu'on laisse $D$ et $E$ mobiles, on est amené à envisager un diorisme suivant la position de $D$. Quand $D$ est situé dans un petit voisinage de $A$, on va utiliser les propriétés de la tangente en $A$ pour étudier la corde $DA$. Pour les autres positions de $D$, il faudra trouver un autre raisonnement.

    Donnons-nous un voisinage de $A$ sur le cercle. Autrement dit, soit $I'$ un point du cercle, point qu'on imagine proche de $A$. Si $D\in\wideparen{AI'}$, on a la majoration grossière $AU>AN$ (car $\wideparen{ANU}$ est obtus). Il reste à montrer que
    $$\frac{\wideparen{CE}}{CJ}\stackrel{?}{>}\frac{\wideparen{AD}}{AN}.$$
    Il suffit désormais de montrer qu'on peut choisir $I'$ suffisamment proche de $A$ pour que cette inégalité soit vérifiée, c'est-à-dire tel que
    $$\frac{CJ}{AN}<\frac{\wideparen{CE}}{\wideparen{AD}},$$
    ou encore, en notant $\dfrac{\alpha}{\beta}=\dfrac{\wideparen{CE}}{\wideparen{AD}}$, tel que~:
    $$AN>\frac{\beta}{\alpha}CJ.$$
    Ibn al-Haytham construit un point $O$ sur la tangente en $A$ tel que
    $$AO=\frac{\beta}{\alpha}CJ,$$
    puis il utilise la majoration grossière $AN>AO$ si $N$ est le projeté de $O$ sur l'axe, parallèlement à $AC$. Il faut alors quand-même vérifier que $O$ est bien situé au-dessus de l'axe, comme on l'a dessiné fig.~8. Pour vérifier ceci, il faut détenir une majoration de $\dfrac{\beta}{\alpha}=\dfrac{\wideparen{AD}}{\wideparen{CE}}$~; pour ce faire, on n'a qu'à utiliser le lemme 1. Nous renvoyons au texte original pour les détails~: cf. \textit{Math. inf.} vol.~V p.~744, où Ibn al-Haytham obtient 
    $$\frac{\beta}{\alpha}=\frac{AT}{TC}\in\left[\frac{AH}{HC},\frac{AS}{SC}\right],$$
    et où il montre ainsi que $T$ est entre $H$ et $S$, et donc que $O$ est au-dessus de l'axe (voir fig.~3 p.~744 pour les notations). La première partie du diorisme est donc démontrée, et on peut d'ores et déjà conclure que~:
    $$(\forall y\in\;]0,\pi/2[)\ (\exists\eta>0)\ (\forall x\in\;]y-\eta,y[)\ f(x)>f(y).$$
    On dira que, pour tout $y\in\;]0,\pi/2[$, la fonction $f$ est \emph{localement décroissante à gauche} en $y$.
    
    Hélas, la seconde partie du diorisme est moins convaincante dans le traité d'Ibn al-Haytham. Sa démonstration consiste, semble-t-il, à utiliser ce résultat local, conjointement au cas particulier du lemme 2, pour obtenir un énoncé un peu plus général, de la manière suivante. Si $0<x<y<z\leq\pi/2$, $f(x)>f(y)$ et $f(y)>f(z)$, alors $f(x)>f(z)$. Or Ibn al-Haytham a montré (lemme 2)~:
    $$(\forall x<y\leq\pi/2)\quad\frac{\sin x}{\sin(x/2)}>\frac{\sin y}{\sin(y/2)},$$
    d'où en posant $x=cy$ ($c\in\;]0,1[$)~:
    $$(\forall y\leq\pi/2)\quad\frac{\sin(y/2)}{\sin(cy/2)}>\frac{\sin y}{\sin cy}.$$
    Alors si $x\in\left]\dfrac{y-\eta}{2^n},\dfrac{y}{2^n}\right[$, on aura $f(x)>f(2x)>\dots>f(2^nx)>f(y)$. Ainsi~:
    $$\left(\forall x\in\bigcup_{n=0}^\infty\left]\frac{y-\eta}{2^n},\frac{y}{2^n}\right[\right)\quad f(x)>f(y).$$
    Mais si $\eta$ est trop petit, cette réunion d'intervalles ne recouvrira pas tout l'intervalle $]0,y[$. Il restera des lacunes à combler~!

    Cette idée nous en inspire une autre, à laquelle peut-être Ibn al-Haytham avait aussi songé. On pourrait construire, par récurrence, une suite de points $y_n$ en prenant $y_0=y$,\enskip $y_1\in\;]y-\eta,y[$, \dots\ , $y_{n+1}\in\;]y_n-\eta_n,y_n[$, \dots\ , de sorte que~:
    $$(\forall x\in\;]y_n-\eta_n,y_n[)\ f(x)>f(y_n)>f(y_{n-1})>\dots>f(y_1)>f(y).$$
    L'union $\bigcup\limits_{n=0}^\infty]y_n-\eta_n,y_n[$ est alors bien un \emph{intervalle}, mais ici encore, rien n'indique qu'elle recouvre tout l'intervalle $]0,y[$. Il est possible que la suite $(y_n)$ ait un point d'accumulation $y'>0$. On pourrait certes recommencer en partant de $y'$, construire une autre suite $(y'_n)$, et ainsi de suite\dots\ On heurte ici un problème d'induction transfinie~: sans un axiome garantissant la possibilité d'une telle induction transfinie sur la droite réelle, rien ne permettait alors de justifier un tel procédé.

    Aujourd'hui, on peut conclure au moyen de la théorie du continu développée au dix-neuvième siècle en démontrant la proposition suivante~: \emph{si $f$ est localement décroissante à gauche en tout point\,\footnote{Au sens défini ci-dessus, c'est-à-dire que, pour une fonction $f$ définie sur un intervalle $I$ ouvert~: $(\forall y\in I)\ (\exists \eta>0)\ (\forall x\in\;]y-\eta,y[)\ f(x)>f(y)$.}, et continue, alors $f$ est décroissante.} Démontrons cette proposition. On fixe $y$. La fonction $f$ est définie sur un intervalle $I$ contenant $y$. On pose $x_0=\sup\lbrace x\in I\mid x<y\text{ \enskip et\enskip}f(x)<f(y)\rbrace$. Supposons que $x_0\in I-\partial I$. Puisque $f$ est localement décroissante à gauche en $y$, il est clair que $x_0<y$, et l'existence du $\sup$ est une propriété de $\mathbb{R}$. La continuité de $f$ à gauche en $x_0$ montre que $f(x_0)\leq f(y)$. La définition du $\sup$ implique que $(\forall x>x_0)\ f(x)\geq f(y)$, et la continuité de $f$ à droite en $x_0$ permet d'en déduire que $f(x_0)\geq f(y)$. Alors $f(x_0)=f(y)$. Mais $f$ est localement décroissante à gauche en $x_0$~:
    $$(\exists\eta>0)\ (\forall x\in]x_0-\eta,x_0[)\ f(x)>f(x_0)=f(y).$$
    Il faut donc que $\sup\lbrace x\in I\mid x<y\text{ \enskip et\enskip}f(x)<f(y)\rbrace\leq x_0-\eta$. Absurde~; donc $x_0\not\in I-\partial I$, ce qu'il fallait démontrer.

    Le lecteur impatient n'aura pas tort s'il réplique~: à tout prendre, $f(x)=\dfrac{\sin x}{\sin cx}$ est continue, mais elle est même dérivable~! Si $f$ est dérivable et localement décroissante en tout point, alors sa dérivée est négative en tout point, et le théorème des accroissements finis implique qu'elle est décroissante. Bien sûr ce théorème dépend du théorème de la valeur intermédiaire, donc aussi d'une propriété forte concernant la nature transfinie du continu, outre le concept même de dérivée, toutes choses qui pourraient sembler déplacées dans le contexte des mathématiques médiévales\dots\ Et pourtant, s'il était impossible qu'Ibn al-Haytham ait recours au concept de dérivée qui n'existait pas encore, il n'en est pas de même du concept de continuité qui pouvait au moins être pensé par le truchement du \emph{mouvement} continu.

    Ibn al-Haytham ne se contente pas d'affirmer l'existence d'un $\eta$, c'est-à-dire du point $I'$ sur la fig.~8. Il construit $I'$ par une propriété d'incidence. $I'$ dépend entièrement de $A$ et de $C$ par une construction à la règle et au compas. On peut donc imaginer une machine (\textit{organum}) entraînant $I'$ le long de l'arc de cercle quand $A$ se déplace. Je fixe $D$~; je pars de la position $A=D$. Alors $I'$ est à gauche de $D$. Je déplace $A$ vers la droite. $I'$ va rester quelque temps à gauche de $D$, \emph{par continuité du mouvement}. La fonction $f$ est donc aussi localement décroissante \emph{à droite}. De plus, je peux parcourir tout arc $[x,y]\subset\,]0,\pi/2[$ avec le point $A$. J'ai donc un recouvrement ouvert de l'arc de cercle $[x,y]$ par de petits intervalles $]y_k-\eta_k,y_k+\eta_k[$. Je peux en extraire un sous-recouvrement fini. Quitte à réordonner les $y_k$, les intervalles successifs doivent se chevaucher. Pour chaque $k$, on choisit $z_k\in\;]y_k,y_k+\eta_k[\,\cap\,]y_{k+1}-\eta_{k+1},y_{k+1}[$. On obtient ainsi une suite de points intermédiaires entre $x$ et $y$ qui permettent d'écrire~:
        $$f(x)>f(z_0)>f(y_1)>f(z_1)>\dots>f(z_{n-1})>f(y_n)>f(z_n)>f(y),$$
        Donc $f$ est bien décroissante sur $]0,\pi/2[$.

        Si nous avons décrit ce raisonnement certes impossible à écrire dans le langage de la géométrie euclidienne du onzième siècle, c'est pour montrer que le lemme d'Ibn al-Haytham découle de la locale décroissance à gauche de $f$ d'une manière somme toute assez intuitive, grâce au mouvement, une fois admise la continuité de $f$. Il reste toutefois un défaut dans le raisonnement ci-dessus. Le point $I'$ ne dépend pas seulement de $A$, il dépend aussi du point $C$ tel que $\dfrac{\wideparen{AB}}{\wideparen{BC}}=\dfrac{1}{c}$. Nous allons montrer qu'une construction géométrique du point $C$ en fonction du point $A$ ne peut exister que si ce rapport est rationnel.

        En effet, s'il existe une construction géométrique de point $C$ en fonction du point $A$, on doit avoir une relation algébrique de la forme
$${(\forall x\in\mathbb{R})\ R(\sin x,\sin cx)=0}.$$
Alors $R(\sin(x+2k\pi),\sin(c(x+2k\pi)))=R(\sin x,\sin(c(x+2k\pi)))=0$, pour tout entier $k$~; mais le polynôme $R(\sin x,\cdot)$ ne peut avoir qu'un nombre fini de racines, et il faut donc que $c(x+2k\pi)$ ne prenne qu'un nombre fini de valeurs modulo $2\pi$, donc aussi $2ck\pi$. Alors $c\in\mathbb{Q}$.

        Si donc $c\not\in\mathbb{Q}$, alors il n'existe aucun \textit{organum} engendrant le mouvement du point $I'$ à partir du mouvement du point $A$. En revanche, on peut facilement imaginer une construction à la règle et au compas pour tout $c$ de la forme $c=m/2^n\in\mathbb{Q}$, donc prouver la continuité de $f$ et la vérité du lemme pour toutes ces valeurs de $c$ par le raisonnement ci-dessus. Il resterait à généraliser pour $c$ quelconque par une démonstration apagogique\,\footnote{Nous n'avons pas cherché à rédiger une telle démonstration apagogique~; mais Ibn al-Haytham a su, ailleurs, user de démonstrations apagogiques pour produire des énoncés universels sur des inégalités entre rapports, ainsi dans \textit{La Configuration des mouvements}, prop.~6 p.~310-313.}.
        
        \subsection{Démonstration globale}
Ibn al-Haytham a trouvé un autre recours pour démontrer ce lemme dans son traité --~sûrement ultérieur~-- sur \textit{La configuration des mouvements}. Il y démontre le lemme par des encadrements judicieux, de manière globale, sans faire appel à aucune propriété différentielle ou locale. Sa démonstration ne repose que sur la loi des sinus dans les triangles --~elle-même équivalente à nos formules d'addition pour les fonctions sinus et cosinus~-- et sur une majoration d'aire semblable à celle utilisée par Ptolémée pour démontrer le cas limite. On va récrire la démonstration d'Ibn al-Haytham de manière purement algébrique\,\footnote{Nous nous sommes abstenus de trop simplifier cette rédaction, pour ne pas obscurcir la parenté formelle déjà ténue avec le texte d'Ibn al-Haytham (\textit{Math. inf.} vol.~V p.~272-279).}, pour mettre en lumière ces seules relations de dépendance logique. \`A cet effet, nous aurons besoin du lemme trigonométrique suivant, dont nous omettons la démonstration\,\footnote{Ce lemme trigonométrique se démontre de manière purement calculatoire en appliquant les formules d'addition, pour tous $p,q,a,b\in\mathbb{R}$.}~:
        $$\text{Si }\frac{\sin(a+q)}{\sin(b-q)}=\frac{\sin(a-p)}{\sin(b-p)},\text{ alors }\frac{\sin(p+q)}{\sin(p-q)}=1-\frac{\sin(b-a)\sin(p+\pi/2)}{\sin(b-p)\sin(a+\pi/2)}.$$

        Soit donc $0<x<y<\pi/2$. On souhaite montrer que
        $$\frac{\sin x}{\sin cx}\stackrel{?}{>}\frac{\sin y}{\sin cy}.$$
        Posons $z<\pi/2$ tel que
        $$\frac{\sin y}{\sin cy}=\frac{\sin x}{\sin z}.$$
        On a donc $z<x$ et $z<cy$. Il s'agit de démontrer que $z>cx$. Or~:
        $$z>cx\ \iff\ \frac{z}{cy}>\frac{x-z}{(1-c)y}$$
        On applique le lemme trigonométrique ci-dessus aux données suivantes~:
        $$p=\frac{1}{2}((1-c)y+x-z),\quad q=\frac{1}{2}(-(1-c)y+x-z),$$
        $$a=\frac{1}{2}((1-c)y+x+z),\quad b=cy+\frac{1}{2}((1-c)y+x-z).$$
        On obtient~:
        $$\frac{\sin(x-z)}{\sin((1-c)y)}=1-\frac{\sin(cy-z)\sin\left(\frac{1}{2}((1-c)y+x-z)+\frac{\pi}{2}\right)}{\sin(cy)\sin\left(\frac{1}{2}((1-c)y+x+z)+\frac{\pi}{2}\right)}.$$
        On vérifie que $(1-c)y+x-z\leq x+y\leq\pi$\enskip et\enskip $(1-c)y+x+z\leq x+y\leq\pi$, de sorte que tous les sinus dans cette formule sont positifs. En particulier, on en déduit que
        $$\frac{\sin(x-z)}{\sin((1-c)y)}<1,$$
        d'où $x-z<(1-c)y$. On peut alors appliquer le cas limite de Ptolémée~:
        $$\frac{x-z}{(1-c)y}<\frac{\sin(x-z)}{\sin((1-c)y)}.$$
        De plus,
        $$\frac{\pi}{2}<\frac{1}{2}((1-c)y+x-z)+\frac{\pi}{2}<\frac{1}{2}((1-c)y+x+z)+\frac{\pi}{2}<\pi.$$
        d'où, dans l'égalité ci-dessus,
        $$\frac{\sin(x-z)}{\sin((1-c)y)}<1-\frac{\sin(cy-z)}{\sin(cy)}.$$
        Appliquons à nouveau le cas limite de Ptolémée aux arcs $cy-z<cy$~:
        $$\frac{\sin(cy-z)}{\sin(cy)}>\frac{cy-z}{cy}.$$
        On a donc~:
        $$\frac{x-z}{1-cy}<1-\frac{\sin(cy-z)}{\sin(cy)}<1-\frac{cy-z}{cy}=\frac{z}{cy}.$$
        C'était la majoration souhaitée.
        
\end{document}